\documentclass[onefignum,onetabnum]{siamart190516}
\pdfoutput=1

\usepackage{lipsum}
\usepackage{amsfonts}
\usepackage{graphicx}
\usepackage[caption=false,farskip=0pt]{subfig}
\usepackage{epstopdf}
\usepackage{algorithmic}
\ifpdf
  \DeclareGraphicsExtensions{.eps,.pdf,.png,.jpg}
\else
  \DeclareGraphicsExtensions{.eps}
\fi

\newsiamremark{remark}{Remark}
\newsiamremark{hypothesis}{Hypothesis}
\crefname{hypothesis}{Hypothesis}{Hypotheses}
\newsiamthm{claim}{Claim}
\def\dfrac{\displaystyle\frac}

\headers{Barycentric Interpolation Based on Equilibrium Potential}{Kelong Zhao, and Shuhuang Xiang}

\title{Barycentrical Interpolation Based on Equilibrium Logarithmic Potential\thanks{Submitted to the editors DATE.
\funding{This work was funded by National Science Foundation of China (No. 12271528).}}}

\author{Kelong Zhao\thanks{School of Mathematics and Statistics, Central South University, Changsha 410083, Hunan, People's Republic of China 
  (\email{clonezhao.1994@gmail.com}).}
\and Shuhuang Xiang\thanks{School of Mathematics and Statistics, Central South University, Changsha 410083, Hunan, People's Republic of China 
  (\email{xiangsh@csu.edu.cn}) (corresponding author).}}

\usepackage{amsopn}

\ifpdf
\hypersetup{
  pdftitle={Barycentric Interpolation Based on Equilibrium Logarithmic Potential},
  pdfauthor={Kelong Zhao, and Shuhuang Xiang}
}
\fi

\begin{document}

\maketitle

\begin{abstract}
  We present a novel barycentric interpolation algorithm designed for analytic functions $f\in\mathcal{A}(E)$ defined on the complex plane. The algorithm, which encompasses both polynomial and rational interpolation, is tailored to handle singularities near $E$. Our method is applicable to regions $E$ bounded by piecewise smooth Jordan curves, and it imposes no connectivity restrictions on the region. The key feature of our approach lies in efficiently computing discrete points via the numerical solution of Symm's integral equation, enabling the construction of polynomial or rational barycentric interpolants. Furthermore, our method provides relevant parameters for the equilibrium potential, such as Robin's constant, which can be used to estimate convergence rates. Numerical experiments demonstrate the convergence rate achieved by our method in comparison to the theoretical convergence rate.
\end{abstract}

\begin{keywords}
  barycentric interpolation, rational approximation, polynomial approximation, potential theory, equilibrium potential
\end{keywords}

\begin{AMS}
  30C10, 30E10, 41A20, 65E05
\end{AMS}
\section{Introduction}
\label{sec:int}
Interpolation stands out as one of the most fundamental and widely employed techniques in scientific computing. Polynomial and rational interpolation, in particular, serve as cornerstone approximation methods across a broad spectrum of numerical analysis \cite{Berrut2014,Berrut2004,Floater2007,Trefethen2013}. However, the traditional scope of these methods is often confined to intervals or circles, presenting challenges in their extension to arbitrary complex plane regions. Recent years have witnessed a growing interest in rational approximation within the broader context of the complex plane, along with its associated challenges and applications \cite{Driscoll2024,Gopal2019,Gopal20192,Nakatsukasa2018}. Therefore, this paper aims to introduce an interpolation algorithm that harnesses the strengths of traditional methods while offering applicability to general complex plane regions.

In this paper, we introduce polynomial and rational interpolants tailored for application on a bounded region $E$ with piecewise smooth boundaries in the complex plane. The polynomial interpolation method is well-suited for both singly connected and disconnected regions, providing an optimal convergence rate for the analytic function $f$. In contrast, rational interpolation is adept at handling singularities near $E$ and can be extended to multiconnected regions provided $f$ exhibits singularities within the holes of the region. Collectively, we refer to these interpolations as Barycentric Interpolation based on Equilibrium Potential (BIEP). Both interpolation techniques are formulated as barycentric formulas and are derived from the logarithmic equilibrium potential.

Logarithmic potential theory has played a crucial role in the development of approximation theory, particularly in the context of polynomial and rational interpolation \cite{Platte2005,TT2010,Trefethen2013}. The logarithmic potential is widely recognized for its ability to characterize the convergence rate of polynomial and rational interpolation for analytic functions \cite{Saff2010,Walsh1965}. While potential theory has traditionally been instrumental in theoretical analysis, the BIEP takes a novel approach by directly computing the discrete point distribution from the density function of the equilibrium potential. This distinctive method sets it apart from existing approaches in the field.

\subsection{Polynomial interpolation of the equilibrium potential}
Let $E$ be a compact set in the complex plane, and let $\mathcal{M}(E)$ denote the collection of all positive unit measures $\mu$ supported on $E$. Then the logarithmic potential with respect to $\mu$ is given by
\[
U_{\mu}(z) = \int \log \frac{1}{|z - t|} \,\mathrm{d}\mu(t) = \int \log \frac{1}{|z - t|} w(t) \,\mathrm{d}t,
\]
where $w$ is the density function of $\mu$. The energy of the logarithmic potential is defined as
\[
I_{\mu}(E) = \iint \log \frac{1}{|t - z|} \, \mathrm{d}\mu(t) \mathrm{d}\mu(z),
\]
and by extension, the Robin constant of $E$ can be defined as
\[
V_E = \inf\{I(\mu):\mu\in\mathcal{M}(E)\}.
\]
The logarithmic capacity of $E$ is defined as
\(\text{cap}(E) = \exp(-V_E)\),
and it is defined as 0 when $V_E = +\infty$.

If $\text{cap}(E) > 0$, there exists a unique measure $\mu_E \in \mathcal{M}(E)$ such that $I(\mu_E) = V_E$. 
This measure is known as the equilibrium measure of $E$. The logarithmic potential corresponding to the equilibrium 
measure is termed the equilibrium potential and satisfies
\(U_{\mu_E}(z) =  V_E \,\, \text{q.e.}\) on $E$, 
where ``q.e." (quasi-everywhere) implies that there exists a set $K$ of capacity 0 such that this property holds 
in $E \setminus K$. It is important to note that if the measure $\mu$ is invariant on $E$, then $\mu = \mu_E$. 
A comprehensive introduction to these fundamental concepts can be found in \cite[Section 1]{Levin&Saff2006}.

Certain discrete point distributions $\{x_i^{(n)}\}_{i=0}^n\subset E$ are considered ``good" and are associated with 
the equilibrium measure $\mu_E$ in a limiting sense. A sequence $\{\mu_n:=\sum_{i=0}^{n}\delta_{x_i^{(n)}}/(n+1)\}$ 
is said to converge weak* to $\mu$ (denoted as $\mu_n\stackrel{*}{\rightarrow}\mu$) if the unit countable 
measure $\mu_n$ satisfies $\int f \, \mathrm{d}\mu_n\rightarrow\int f \, \mathrm{d}\mu$ as $n\rightarrow\infty$
for all continuous functions $f$. For instance, the unit countable measures of Jacobi points on $[-1,1]$ 
converge weak* to the equilibrium measures on $[-1,1]$. The same weak* convergence holds for Fekete points 
as well as for Leja points \cite{Narayan2014}.

\begin{figure}[htbp]
  \centering
  \includegraphics[width=0.99\linewidth]{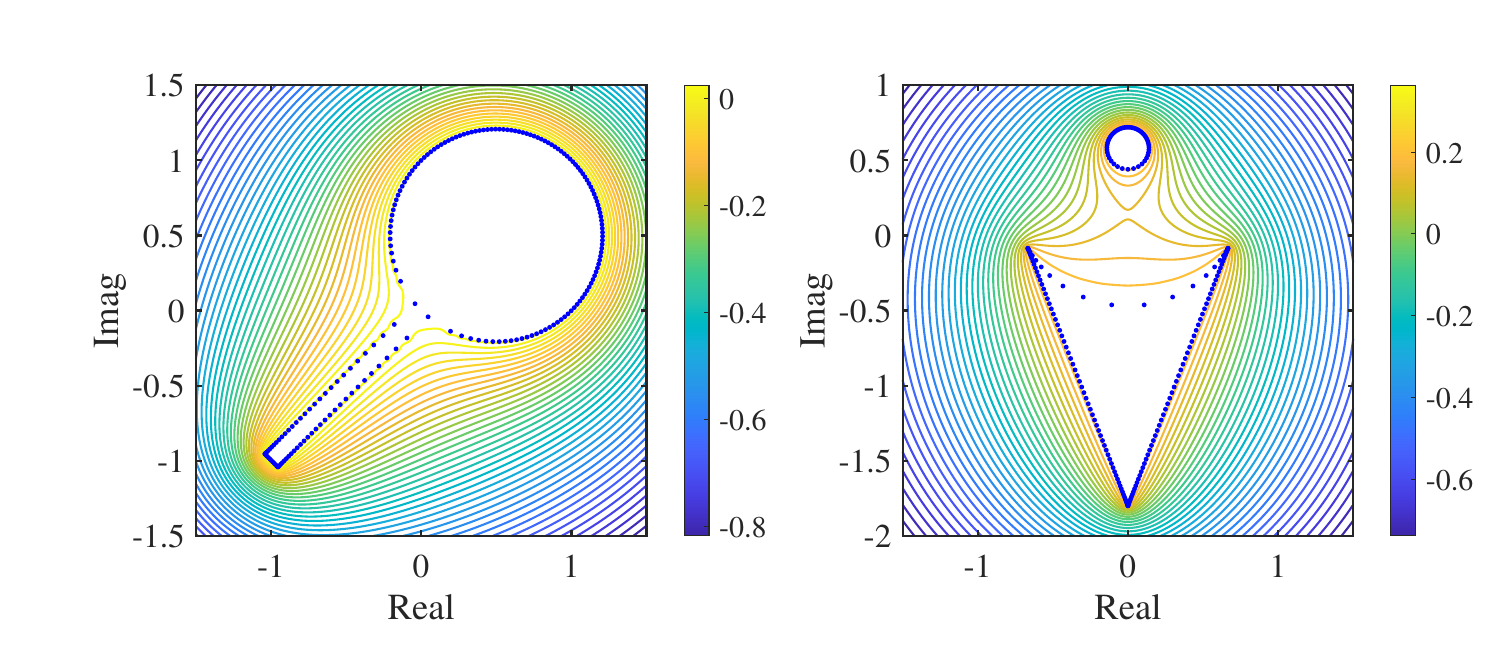}
  \caption{The contours of the potential $U_{\mu_n}$ ($n=200$) in the `lollipop' (left) and `ice cream cone' (right) domains. 
  These points (blue) resemble identically charged particles that repel each other and exhibit tip aggregation.}
  \label{Figure_1}
\end{figure}

Unlike Fekete and Leja points, which require solving a nonlinear optimization problem, we compute the discrete 
point distribution directly from the density function corresponding to the equilibrium measure. For singly connected 
and disconnected regions with piecewise smooth boundaries, the discrete points $\{z_i\}$ can be computed efficiently. 
\cref{Figure_1} illustrates these points and the corresponding contours of the potential $U$. These 
contours surround $E$, but the interior of $E$ is almost empty. This indicates that the potential is nearly 
constant in the interior of $E$, similar to the characteristics of an equilibrium logarithmic potential.

These discrete points of the approximate equilibrium measure can be used directly as interpolation nodes to 
construct polynomial interpolants. When the nodes are given, the polynomial barycentric interpolation \cite{Berrut2004} is given by
\begin{equation} \label{eq:1.1p}
  p_n(x)={\sum\limits_{k=0}^{n}\dfrac{w_k}{x-x_k}f(x_k)}\Big/{\sum\limits_{k=0}^{n}\dfrac{w_k}{x-x_k}}, \,\, w_k=C/\prod_{j=0,j\ne k}^n(x_k-x_j),
\end{equation}
where $C$ is a nonzero constant.
We will detail the computational procedure for this polynomial interpolant in \cref{sec:3.1}.

An approximation of the Robin constant $V_E$ can also be obtained in the polynomial interpolation of BIEP. 
The Robin constant can be used to estimate the convergence rate of the optimal 
polynomial interpolation on a compact set of the complex plane $E$ \cite[Section 2]{Levin&Saff2006}. 
Specifically, if $C\setminus E$ is connected and $f$ is analytic on $E$, then 
\begin{equation}\label{eq:1.1e}
  \limsup_{n\to\infty} \left(\min_{p\in\mathcal{P}_n}\|p-f\|_E \right)^{\frac{1}{n}}=\exp\left\{ U^{\mu_E}(\Gamma_R)-V_E\right\}
\end{equation}
where $\Gamma_R$ is one of the contours of $U^{\mu_E}$ and is the largest contour such that $f$ is analytic within $\Gamma_R$. 
It is worth noting that the numerical examples in \cref{sec:4.1} show that our polynomial interpolant exhibits 
a convergence rate that approximates the optimal rate.

\subsection{Rational interpolation of the equilibrium potential}
When the compact set $E$ is given, Eq. \eqref{eq:1.1e} shows that the rate of polynomial interpolation depends on the value
of the logarithmic potential on the contour $\Gamma$. This value is determined by the location of the singularities of $f$. 
When the value of the logarithmic potential at these singularities is close to the Robin constant $V_E$, the 
convergence of polynomial interpolation may be slow. This type of near-singular function is often seen in problems 
such as singularly perturbed problems \cite{Chen2011,McCoid2019,Tang1996}, finite time blow-up \cite{xia2021,Yang2013}, 
and problems with small radii of curvature or almost-singular boundary data \cite{Gopal2019,Helsing2008,Hochman2013}. 
Rational functions are commonly used to approximate near-singular functions in such problems.

Similar to polynomials, rational interpolation is closely related to the equilibrium potential but belongs to the 
logarithmic potential with signed measures. Let $E, F \subset \mathbb{C}$ be two closed sets that are a positive distance 
apart, and $\mu_E, \mu_F$ be positive unit measures supported on $E$ and $F$, respectively. The logarithmic potential 
on the signed measure $\mu=\mu_E-\mu_F\in \mathcal{M}(E,F) $ can be defined as
\[
U_\mu(z) = \int \log \frac{1}{|z - t|}\, \mathrm{d}\mu(t).
\]
The energy $I$ and minimum energy $V$ can be defined as
\[
I(\mu) = \iint \log \frac{1}{|z - t|}\, \mathrm{d}\mu(z)\, \mathrm{d}\mu(t)\quad \text{and}\quad V(E,F) = \inf_{\mu' \in \mathcal{M}(E,F)} I(\mu'),
\]
respectively. The capacity $\text{cap}(E,F)$ between $E$ and $F$ is given by
\[
\text{cap}(E,F) = 1/V(E,F).
\]
The measure that achieves the minimum energy is referred to as the equilibrium measure $\mu^*$. 
It is worth noting that only the equilibrium potential satisfies
\begin{equation}\label{eq:2g}
U_{\mu^*}(z) = \begin{cases}
  c_1, \quad z\in E\\
  -c_2,\quad z\in F
\end{cases}, 
\end{equation}
where constants $c_1$ and $c_2$ have $c_1+c_2=1/\text{cap}(E,F)$.
A detailed description of these basic concepts and properties can be found in \cite[Section 6]{Levin&Saff2006}.

\begin{figure}[hpbt]
  \centering
  \includegraphics[width=0.55\linewidth]{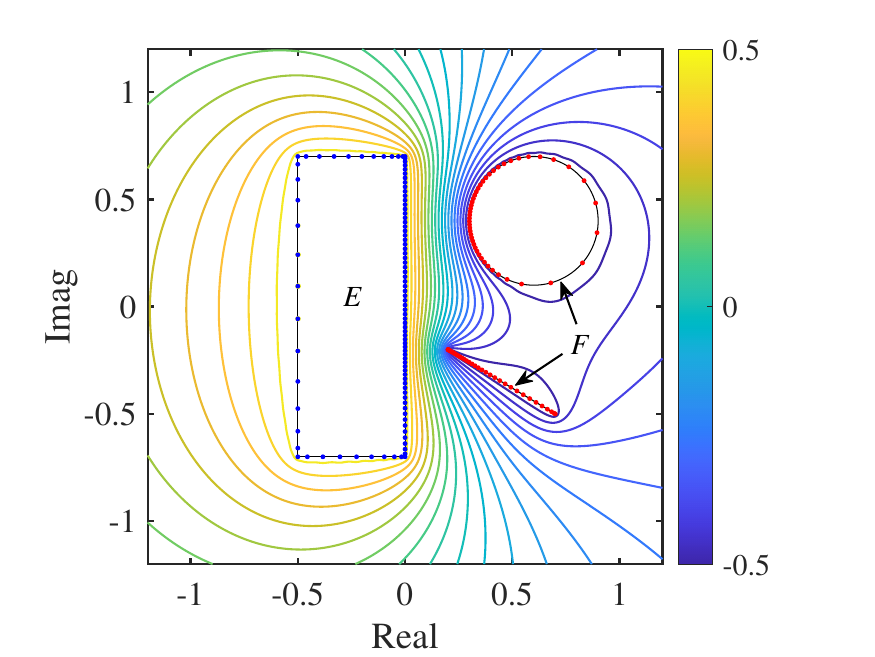}
  \caption{By our method, the nodes (blue) and poles (red) are selected on $\partial E$ and $\partial F$, respectively. These nodes $\{x_i\}_{i=0}^n$ and poles $\{z_j\}_{j=1}^n$ make the discrete potential 
  $U_{\mu_n}$ almost constant on $E$ and $F$. Therefore, the interior of $E$ and $F$ is almost blank in the contour map of the potential $U$.}
  \label{Figure_2}
\end{figure}

Let $\{x_i^{(n)}\}_{i=0}^n\subset E$ be the nodes and $\{z_j^{(n)}\}_{j=1}^n\subset F$ be the poles. We anticipate that the signed measure
\[
\mu_n = \frac{1}{n+1}\sum_{i=0}^n \delta_{x_i^{(n)}} - \frac{1}{n+1}\sum_{j=1}^n \delta_{z_j^{(n)}}
\]
weak*-converges to the equilibrium measure $\mu^*_{E,F}$. \cref{Figure_2} depicts the nodes (blue) and poles (red) 
computed by the method proposed in this paper, along with the contours of the potential $U_{\mu_n}$. 
The contours enclose $E$ and $F$, indicating that the logarithmic potential is approximately invariant on $E$ and $F$, 
consistent with the behavior described in \eqref{eq:2g}.

When the nodes and poles are known, rational barycentric interpolation can be calculated using \eqref{eq:1.1p}, where the barycentric weights are given by
\begin{equation} \label{eq:1.1c}
  w_k = C \dfrac{\prod_{j=1}^m(x_k-z_j)}{\prod_{i=0,i\ne k}^n(x_k-x_i)},\,\, C\neq0\,\, (m\leq n),
\end{equation}
as proposed by Berrut \cite{Berrut1997}. We provide a detailed procedure for rational interpolation of BIEP in \cref{sec:3.2}.

The convergence rate of rational interpolation can also be expressed in terms of a logarithmic potential. Let $R_n$ denote the set of rational functions whose numerators and denominators do not exceed degree $n$. Then, we have
\begin{equation}\label{eq:1.2c}
\limsup_{n\to\infty} \left(\inf_{r\in R_n} \|r-f\|_E \right)^{\frac{1}{n}}\le\exp\left\{ -1/\text{cap}(E,\partial \Omega)\right\},
\end{equation}
if $f$ is analytic on $\Omega$ with $E\subset\Omega$ \cite[Section 6]{Levin&Saff2006}. 
For some special cases where $\partial F$ is close to $\partial\Omega$, the rate of asymptotic convergence can be approximated as $[\exp(-c_1-c_2)]^n$.
In this paper, we will use $[\exp(-c_1-c_2)]^n$ as a reference convergence rate. 
However, the rational interpolation we provide tends to converge slightly slower than this in the case of branch singularities.

We will test the effectiveness of the approximation of near-singular functions and compare it with the theoretical 
convergence rate in \cref{sec:4.2}. \Cref{sec:discu} will discuss the choice of parameters and 
present more complex examples. In \cref{sec:Laplace}, we will demonstrate a simple application of BIEP,
namely solving 2-D Laplace equations. Finally, \cref{sec:conclusion} will discuss 
further potential uses of BIEP.
\section{Algorithms of the interpolation}\label{sec:algorithm}
The connection between polynomial and rational interpolation and the logarithmic equilibrium potential lies 
in the weak* convergence of the discrete measure $\mu_n$ to the equilibrium measure $\mu^*$. However, the 
definition of weak* convergence does not directly apply to the construction of the discrete measure $\mu_n$. 

It is feasible to inscribe the distribution of nodes and poles by a density function on the boundary. An intuition 
is that high-density sections distribute more discrete points and the opposite is true for low-density. We can define 
the relationship between discrete points and the density function in this way: 
\begin{definition}\label{def:1}
A family of point sets $\{\{x_k^{(n)}\}_{k=0}^n:n=1,2,\cdots\}$ obeys the density function $w(t)>0$ for all $t\in\partial E$ if for any segment $\partial E[\widehat{a,b}] \subseteq \partial E$, the family of point sets satisfies
\begin{equation} \label{eq:3a}
\lim\limits_{n\to\infty}\frac{n_{\partial E[\widehat{a,b}]}}{n+1}=\int_{\partial E[\widehat{a,b}]}w(t)\,\lvert\mathrm{d}t\rvert,
\end{equation}
where $a$ and $b$ are the two endpoints of the segment $\partial E[\widehat{a,b}]$ and $n_{\partial E[\widehat{a,b}]}$ denotes the number of the points on $\partial E[\widehat{a,b}]$.
\end{definition}

Consistent with the intuition, the countable measure corresponding to a sequence of discrete points satisfying 
\cref{def:1} converges weak* to the continuous measure corresponding to that density function: 
\begin{theorem}\label{thm:den}
If $\partial E$ is a bounded piecewise simply smooth boundary and a family of point sets $\{x_i^{(n)}\}_{i=0}^n\subseteq \partial E$ obeys a positive density function $w$ of a unit measure $\mu$ on $\partial E$, then it holds
\begin{equation}\label{eq:pdis}
\lim_{n\to\infty}\max_{1\leq i\leq n}\lvert x_{i-1}^{(n)}-x_i^{(n)}\rvert = 0
\end{equation}
and $\mu_n\stackrel{*}{\rightarrow}\mu$.
\end{theorem}

For consistency of the implementation of the rational interpolation, we sketch the proof of \cref{thm:den} in Appendix A.

In this paper, the family of the point sets is chosen such that for any two adjacent points $x_i^{(n)},x_{i+1}^{(n)}\in\{x_k^{(n)}\}_{k=0}^n$, it holds
\[\int_{\partial E[\widehat{x_i^{(n)},x_{i+1}^{(n)}}]}w(t)\,\lvert\mathrm{d}t\rvert= \begin{cases}
	1/(n+1), & \text{if } \partial E \text{ is a closed curve, } i = 0,1,\cdots,n \\
	1/n, & \text{if } \partial E \text{ is a curve segment, } i = 0,\cdots,n-1
\end{cases},\]
where $x_{n+1}^{(n)}=x_0^{(n)}$ if $\partial E$ is a closed curve.
It is easy to verify that $\{x_i^{(n)}\}_{i=0}^n$ satisfies \cref{def:1}.
The above equation can be expressed as
\begin{equation} \label{eq:3b}
	\int_{\partial E[\widehat{x_0^{(n)},x_i^{(n)}}]}w(t)\,\lvert\mathrm{d}t\rvert = \begin{cases}
		i/(n+1), & \text{if } \partial E \text{ is a closed curve} \\
		i/n, & \text{if } \partial E \text{ is a curve segment}
	\end{cases},\ i = 0,1,\cdots,n.
\end{equation}
If $\partial E$ is a closed curve, $x_0^{(n)}$ could be any point on $\partial E$.
If $\partial E$ is a curve segment, $x_0^{(n)}$ is an endpoint of the curve segment.

As an example, we consider the density function \(w(t) = (\pi \sqrt{1-x^2})^{-1}\) on \(E=[-1,1]\). 
Suppose a family of point sets \(\{\{x_k^{(n)}\}_{k=0}^n : n=1,2,\cdots\}\) on \(E\) satisfies \eqref{eq:3b}, 
then it yields \(x_i^{(n)} = \cos\left(\frac{n-i}{n}\pi\right)\), and we obtain the Chebyshev-Lobatto points. 
There exists a rational interpolation that is also consistent with \eqref{eq:3b} \cite{Deun2009}.

However, obtaining an explicit formula for the density function $w$ of the equilibrium measure can be challenging,
especially for general regions $E$ and $F$ with various shapes and relative positions. Moreover, even if an analytical
expression for the density function is obtained, it may not be possible to explicitly express the point distribution
from \eqref{eq:3b}. To overcome these challenges, we will apply a step function $\hat{w}(t)$ to approximate
the density function $w(t)$ by solving a Symm's integral equation.

\subsection{Polynomial interpolation based on equilibrium potential} \label{sec:3.1}
The most straightforward scenario is polynomial interpolation with poles situated at infinity in the complex plane. In this case, we can obtain an approximation of the equilibrium potential's density function $\hat{w}$ and Robin constant $V_E$ by solving Symm's integral equation \cite{Kythe1998}. Symm's equation is a Fredholm integral equation of the first kind, which has a logarithmic singular kernel:
\begin{equation} \label{eq:3c}
  \int_{\partial E} \log\frac{1}{\lvert z-t\rvert}w(t)\,\lvert\mathrm{d}t\rvert=V_E, \quad z\in\partial E,
\end{equation}
where $w$ is the density function of the equilibrium measure and satisfies
\begin{equation} \label{eq:density}
  \int_{\partial E} w(t)\,\lvert\mathrm{d}t\rvert=1.
\end{equation}
Due to the logarithmic kernel being singular, the ill-conditioning of \eqref{eq:3c} is quite manageable \cite[Subsection 1.1.5]{Atkinson1997}.
Numerous numerical methods have been developed to solve Symm's integral equations \cite[Chapter 9]{Kythe1998}. In this work, we utilize the constant element method to solve \eqref{eq:3c} and \eqref{eq:density}. This method is straightforward to implement and provides an approximation in the form of a step function.

We partition $\partial E$ into $N$ segments, $\partial E_1,\cdots,\partial E_N$, using $t_0,\cdots,t_{N}$ as the dividing points. If the curve is closed, the first and last points coincide. It is noteworthy that $N$ does not depend on $n$, but increasing $N$ enhances the precision of the density function approximation, yielding more precise discrete points.

Let $\hat{w}_i(\approx w(t))$ for $t\in \partial E_i$. Taking a point $t_{i+1/2}\in \partial E_i$ in each
subinterval $\partial E_j$, (\ref{eq:3c}) is simplified to
\begin{equation} \label{eq:3d}
  \sum_{j=1}^{N}\hat{w}_j\int_{\partial E_j}\log\frac{1}{\lvert t_{i+1/2}-t\rvert}\,\lvert\mathrm{d}t\rvert=V_E,\  i=1,\cdots,N.
\end{equation}
Then (\ref{eq:3d}) can be represented as $\mathbf{A}\mathbf{\hat{w}}=\mathbf{V_E}$ where
$\mathbf{\hat{w}}=[\hat{w}_1;\cdots;\hat{w}_N]$ and the $a_{i,j}$ of $\mathbf{A}$ denotes
\begin{equation}\label{eq:linapp}
  a_{i,j}\approx\int_{\partial E_j}\log\frac{1}{\lvert t_{i+1/2}-t\rvert}\, \lvert\mathrm{d}t\rvert,\ i=1,\cdots,N,\ j=1,\cdots,N.
\end{equation}
For $i\ne j$, $a_{i,j}$ is evaluated by
\begin{equation*}
  a_{i,j}=-\frac{\lvert\partial E_i\rvert _L}{6}(\log\lvert t_{i+\frac{1}{2}}-t_{j}\rvert
  +4\log\lvert t_{i+\frac{1}{2}}-t_{j+\frac{1}{2}}\rvert+\log\lvert t_{i+\frac{1}{2}}-t_{j+1}\rvert),
\end{equation*}
while for $i=j$,
\begin{equation*}
	a_{i,i}=\lvert t_{i+\frac{1}{2}}-t_{i}\rvert(1-\log\lvert t_{i+\frac{1}{2}}-t_{i}\rvert)
  +\lvert t_{i+\frac{1}{2}}-t_{i+1}\rvert(1-\log\lvert t_{i+\frac{1}{2}}-t_{i+1}\rvert)
\end{equation*}
(see  \cite[P. 241]{Kythe1998}), where $\lvert\partial E_i\rvert _L$ denotes the length of $\partial E_i$.
Then equations \eqref{eq:3c} and \eqref{eq:density} are jointly expressed as
\begin{equation} \label{eq:3e}
  \left( \begin{array}{c|c}
    \mathbf{A} & \mathbf{-1}\\ \hline
    \mathbf{L_E} & 0 \\
  \end{array} \right)
  \left( \begin{array}{c}
    \mathbf{\hat{w}}\\ \hline
    V_E\\
  \end{array} \right) =
  \left( \begin{array}{c}
    \mathbf{0}\\ \hline
    1\\
  \end{array} \right),
\end{equation}
where $\mathbf{L_E}=[\lvert\partial E_1\rvert _L, \cdots,  \lvert\partial E_N\rvert _L]$.
The approximate density $\mathbf{\hat{w}}$ and the Robin constant $V_E$ can be obtained from \eqref{eq:3e}.

Thereafter, we can obtain the interpolation nodes from (\ref{eq:3b}) using the approximate density function $\hat{w}(t)=\hat{w}_i$ for $t\in\partial E_i$.
\begin{figure}[hpbt]
  \centering
  \includegraphics[width=0.75\linewidth]{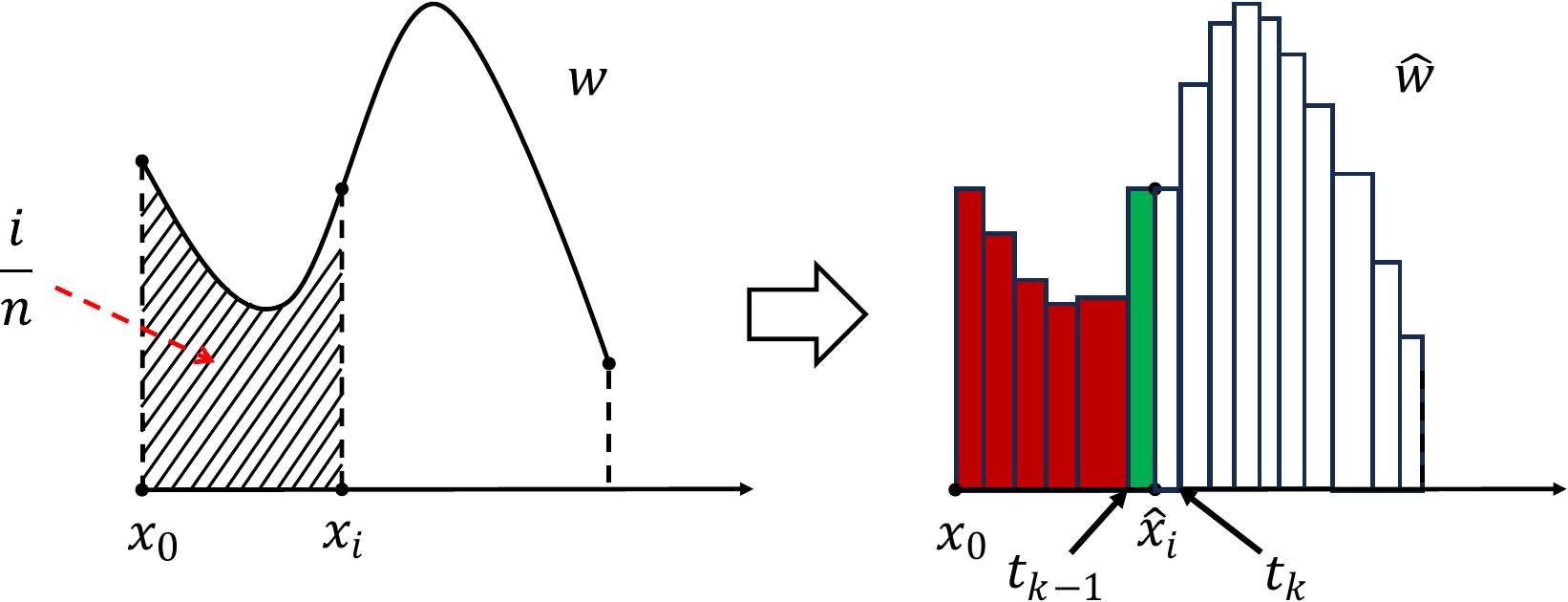}
  \caption{Schematic diagram for $\partial E$ is a curve segment.}
  \label{Figure_3}
\end{figure}

Suppose $\{x_i^{(n)}\}_{i=0}^n\subset \partial E$ and the approximate density function $\hat{w}(t)$ satisfies \eqref{eq:3b}. We set $x_0^{(n)}=t_0$, and define
\[ h_{n,i}=  \begin{cases}
    i/(n+1), & \text{if $\partial E$ is a closed curve}, \\
    i/n, & \text{if $\partial E$ is a curve segment},
  \end{cases} \quad i = 0,1,\cdots,n, \]
and $s_j = \sum\limits_{i=1}^j \lvert\partial E_i\rvert _L\hat{w}_i$ for $j=1,\cdots, N$ $(s_0=0)$. If $s_{k-1}\leq h_{n,i}<s_{k}$, then the approximation $\hat{x}_i^{(n)}\in \partial E[\widehat{t_{k-1},t_{k}}]$, and  \eqref{eq:3b} can be approximated by
\begin{equation}\label{eq:3f}
  \sum\limits_{j=1}^{k-1}\lvert\partial E_j\rvert _L\hat{w}_j+
  \lvert\partial E[\widehat{t_{k-1},\hat{x}_i^{(n)}}]\rvert _L\hat{w}_k = h_{n,i}.
\end{equation}
This allows us to express the length of $\hat{x}_i^{(n)}$ from the point $x_0^{(n)}=t_0$ on $\partial E$ as
\[ \lvert\partial E[\widehat{t_0,\hat{x}_i^{(n)}}]\rvert _L=\sum\limits_{j=1}^{k-1} \lvert\partial E_j\rvert _L + \varepsilon_{i,k},\quad \varepsilon_{i,k} = (h_{n,i}-s_{k-1})/\hat{w}_{k}. \]

Define the matrix $\mathbf{B} = [b_{i,j}]$ for $1\leq i\leq n$ and $1\leq j\leq N$ with
\begin{equation}\label{eq:3g}
  b_{i,j} = \begin{cases}
    1 , & \hat{b}_{i,j} \geq 1 \\
    \hat{b}_{i,j} , & 0 < \hat{b}_{i,j} < 1 \\
    0 , & \hat{b}_{i,j} \leq 0 \\
  \end{cases}, \quad
  \hat{b}_{i,j}=\frac{h_{n,i}-s_{j-1}}{\lvert\partial E_{j}\rvert _L\hat{w}_{j}}.
\end{equation}
We have
\begin{equation}\label{eq:3h}
  [\lvert\partial E[\widehat{t_0,\hat{x}_1^{(n)}}]\rvert _L,\cdots,\lvert\partial E[\widehat{t_0,\hat{x}_n^{(n)}}]\rvert _L]^T
  = \mathbf{B}[\lvert\partial E_1\rvert _L,\cdots,\lvert\partial E_N\rvert _L]^T.
\end{equation}

Given the reference point $t_0$ as $\hat{x}_0^{(n)}$, the nodes $\{\hat{x}_i^{(n)}\}_{i=0}^n$ can be obtained from the curve lengths. Then, we can perform the barycentric interpolation using \eqref{eq:1.1p}.

\begin{algorithm}[H]
  \caption{Barycentric polynomial interpolation with equilibrium potential}\label{alg:A}
  \begin{algorithmic}[1]
    \STATE{$\mathbf{function}$ BPIEP($f \in \mathcal{A}(E)$, $n\in \mathbb{N}$)}
    \STATE{\quad Discretize $\partial E$;}
    \STATE{\quad Get approximate density function $\hat{w}$ by \eqref{eq:3e};}
    \STATE{\quad Compute the nodes $X=\{x_i^{(n)}\}_{i=0}^n$ using Algorithm \ref{alg:B};}
    \STATE{\quad Return $P_n[f]$ from the polynomial case of \eqref{eq:1.1p} with nodes $X$.}
    \STATE{$\mathbf{end\ function}$}
  \end{algorithmic}
\end{algorithm}

\begin{algorithm}[H]
  \caption{Generate points from density function $\hat{w}$}\label{alg:B}
  \begin{algorithmic}[1]
  \STATE{$\mathbf{function}$ den2pts($\hat{w}$,$n\in \mathbb{N}$)}
  \STATE{\quad Normalize $\hat{w}$ and define column vectors $L = [\lvert\partial E[t_i,t_{i+1}]\rvert _L]$, $W = [\hat{w}_i]$;}
  \STATE{\quad Set $S = [0;\lvert\partial E_1\rvert _L\hat{w}_1;\cdots ;\sum_{i=1}^{N-1} \lvert\partial E_i\rvert _L\hat{w}_i]$;}
  \STATE{\quad Compute matrix $\mathbf{B}$ by \eqref{eq:3g};}
  \STATE{\quad Return points $\{\hat{x}_i^{(n)}\}_{i=0}^n$ using \eqref{eq:3h};}
  \STATE{$\mathbf{end\ function}$}
  \end{algorithmic}
\end{algorithm}

It is worth noting that when $E$ is a union of simply connected regions, we tend to divide the boundary of $E$ into multiple parts, where the ratio of the number of points in each part approximates the ratio of the integrals of the density function $\hat{w}$.

\subsection{Rational interpolation based on equilibrium potential} \label{sec:3.2}
To approximate an analytic function $f\in\mathcal{A}(E)$ with singularities near $E$, rational interpolation 
often outperforms polynomial interpolation. The key advantage lies in rational interpolation's ability to introduce 
poles near the singularities of $f$. These additional poles can significantly improve the convergence rate compared to polynomial interpolation methods.

We introduce a region $F$ for distributing the poles of the rational interpolation \eqref{eq:1.1c} to approximate analytic functions $f\in\mathcal{A}(E)$ with singularities near $E$. For specific singularities, we can use simple strategies for selecting $F$.

\begin{itemize}

\item When the singular point $Z$ is isolated, such as a pole or essential singularity, we cover it with a small disk $F=D(Z,r)$ to avoid intersection with $E$ and ensure a lower potential on $F$.

\item For branch singularities, if the two endpoints of a branch cut are $Z_1$ and $Z_2$, and $Z_1$ is close to $E$ but $Z_2$ is far from $E$, we make a line segment from $Z_1$ in the direction away from $E$. We find that a length between $4$ and $10$ is often acceptable when the diameter of $E$ is between $1$ and $2$. For example, set $f(z)=\sqrt{z+0.1}$ and $E=[0,1]$. The branch cut is $[Z_2,Z_1]=[-\infty,-0.1]$ and we select $F=[-4.1,-0.1]$.
If the two ends of the branch cut $[Z_1, Z_2]$ are both near $E$, we typically use the line segment $[Z_1, Z_2]$ as $F$ if it does not intersect $E$. For example, let $f(z) = [(z^2-0.25)/(z^2-0.01)]^{0.5}$, which has four branch points and two branch cuts. Then we choose $F=[-0.5, -0.1]\cup[0.1, 0.5]$. 
In some cases, it may be useful to use a curved segment as $F$ to avoid intersection with $E$. 
\end{itemize}

Similar to the previous subsection, we still start with the potential.
From the properties of the equilibrium potential
and (\ref{eq:2g}), we get
\begin{equation} \label{eq:3i}
  \int_{\partial E}\log\frac{1}{\left\lvert z-t\right\rvert}w_E(t)\,\lvert\mathrm{d}t\rvert-
  \int_{\partial F}\log\frac{1}{\left\lvert z-t\right\rvert}w_F(t)\,\lvert\mathrm{d}t\rvert= \begin{cases}
   c_1, & z\in \partial E\\
   -c_2, & z\in \partial F\\
 \end{cases},
\end{equation}
and
\begin{equation}\label{eq:3j}
\int_{\partial E} w_E(t)\,\lvert\mathrm{d}t\rvert=1,\quad \int_{\partial F} w_F(t)\,\lvert\mathrm{d}t\rvert=1,
\end{equation}
where $w_{S}(t)> 0$, and $S$ denotes $E$ or $F$.
For the more general case, $\partial E$ and $\partial F$ may consist of multiple simple boundaries.

Let $\partial E$ be divided into $I$ subintervals $\{\partial E_k\}_{k=1}^{I}$ and $\partial F$ be divided
into $J$ subintervals $\{\partial F_k\}_{k=1}^J$. 
Similar to the polynomial case, let $\hat{w}_{E_k}(\approx w_E(t))$ for $t\in\partial E_k$ and $\hat{w}_{F_k}(\approx w_F(t))$ for $t\in\partial F_k$,
then \eqref{eq:3i} and \eqref{eq:3j} are jointly expressed as
\begin{equation} \label{eq:3k}
  \left( \begin{array}{c|c|cc}
    \mathbf{A_{E,E}} & \mathbf{-A_{E,F}} & \mathbf{-1} & \mathbf{0} \\ \hline
    \mathbf{A_{F,E}} & \mathbf{-A_{F,F}} & \mathbf{0}  & \mathbf{1} \\ \hline
    \mathbf{L_E}    & \mathbf{0}       & 0           & 0           \\
    \mathbf{0}      & \mathbf{L_F}     & 0           & 0           \\
  \end{array} \right)
  \left( \begin{array}{c}
    \mathbf{\hat{w}_{E}}\\ \hline
    \mathbf{\hat{w}_{F}}\\ \hline
    c_1\\
    c_2\\
  \end{array} \right) =
  \left( \begin{array}{c}
    \mathbf{0}\\ \hline
    \mathbf{0}\\ \hline
    1\\
    1\\
  \end{array} \right)
\end{equation}
where $\mathbf{L_E}=[\lvert\partial E_1\rvert _L, \cdots,  \lvert\partial E_I\rvert _L]$ and $\mathbf{L_F}=[\lvert\partial F_1\rvert _L, \cdots,  \lvert\partial F_J\rvert _L]$.
The element of the $i$-th row and $j$-th column of $A_{E,E}$ $(I \times I)$ is the approximation of the potential generated by $\partial E_j$ under uniform measure at $t_{i+1/2} \in \partial E_i$, as shown in \eqref{eq:linapp}. The element of row $i$ and column $j$ of $A_{E,F}$ $(I \times J)$ is then the approximation of the potential generated by $\partial F_j$ at $t_{i+1/2} \in \partial E_i$ under the uniform measure. 
Similarly, the elements in $A_{F,E}$ $(J \times I)$ and $A_{F,F}$ $(J \times J)$ follow the same rule.

Therefore, we can derive the approximate density function $\hat{w}_S$ used to generate the discrete points.
The parameters $c_1$ and $c_2$ can be used to estimate a theoretical rate of convergence, although not necessarily accurate for the branch singularities.

Using the \texttt{den2pts} (\cref{alg:B}), we can obtain the nodes $\{x_i^{(n)}\}_{i=0}^n$ and poles $\{z_j^{(n)}\}_{j=1}^n$ from the density
function $\hat{w}$. Subsequently, we can obtain the rational interpolation using \eqref{eq:1.1c}. We summarize the rational interpolation with equilibrium potential in \cref{alg:C}.
\begin{algorithm}
  \caption{Barycentric rational interpolation with equilibrium potential}
  \label{alg:C}
  \begin{algorithmic}[1]
  \STATE{$\mathbf{function}$ BRIEP($f \in\mathcal{A}(E)$, singularities $[q_i]$ of $f$, $n\in\mathbb{N}$)}
  \STATE{\quad Select the regions $F$ to cover $[q_i]$;}
  \STATE{\quad Discrete $\partial E$ and $\partial F$;}
  \STATE{\quad Get approximately density functions $\hat{w}_E$ and $\hat{w}_F$ by (\ref{eq:3k});}
  \STATE{\quad The nodes $X=\{x_i^{(n)}\}_{i=0}^n\leftarrow$den2pts($\hat{w}_E$,$n$);}
  \STATE{\quad The poles $Z=\{z_j^{(n)}\}_{j=1}^n\leftarrow$den2pts($\hat{w}_F$,$n-1$);}
  \STATE{\quad Barycentric weights $\{w_k\}_{i=0}^n$ from (\ref{eq:1.1c}) with $X$ and $Z$;}
  \STATE{\quad Return $r_{n}[f]$ with weights $\{w_k\}_{i=0}^n$ and $\{f_i\}_{i=0}^n$.}
  \STATE{$\mathbf{end\ function}$}
  \end{algorithmic}
\end{algorithm}

\section{Experimental results}\label{sec:experiments}
In this section, we examine the impact of polynomial and rational interpolation for BIEP. 
The target functions are analytic, featuring singularities in the vicinity of the region. 
Our focus lies on the approximation error as a function of the order $n$, comparing it against the theoretical 
convergence rate discussed in \cref{sec:int}.

\subsection{Polynomial interpolants}
\label{sec:4.1}
\begin{figure}[htbp]
  \centering
  \includegraphics[width=1\linewidth]{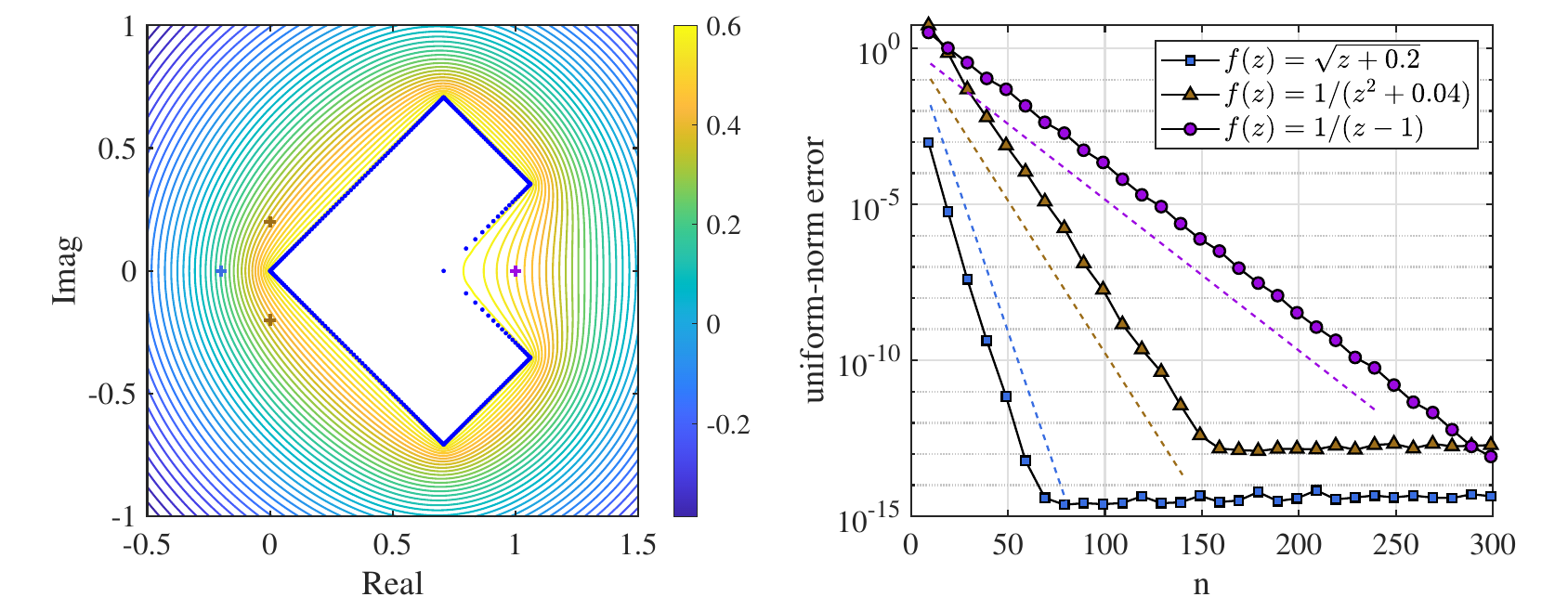}
  \caption{Left: Interpolation nodes (blue) and the discrete potential $U_{\mu_n}$ generated by the 300 nodes. Right:
  The convergence rates compared with the theoretical convergence rates $\rho^n$,
  $\rho=\exp(-0.4180)$ (blue dashed), $\rho=\exp(-0.2248)$ (yellow-brown dashed) and $\rho=\exp(-0.1115)$ (violet dashed), respectively.}
  \label{fig:poly}
\end{figure}

The first example is depicted in \cref{fig:poly}, interpolated over a polygonal region $E$. The nodes are 
distributed so that their potential function is nearly constant over $E$. As $n$ increases, the value of the discrete 
potential function $U_n$ inside $E$ approaches that of $V_E$ computed by \cref{alg:A}. The left panel labels 
the positions of the 300 nodes (blue dots) as well as the positions of the singularities of the three objective 
functions ($+$). The value of the discrete potential $U_n$ is approximately $0.1937$ at $-0.2$ (blue $+$), $0.3868$ 
at $\pm0.2i$ (brown $+$), and $0.5002$ at $1$ (violet $+$). Based on these values, one can estimate the theoretical 
rate of convergence of the three objective functions. The figure on the right illustrates the uniform norm error 
for each of the three objective functions, all of which are consistent with the theoretical rate of convergence. 

\begin{figure}[htbp]
  \centering
  \includegraphics[width=1\linewidth]{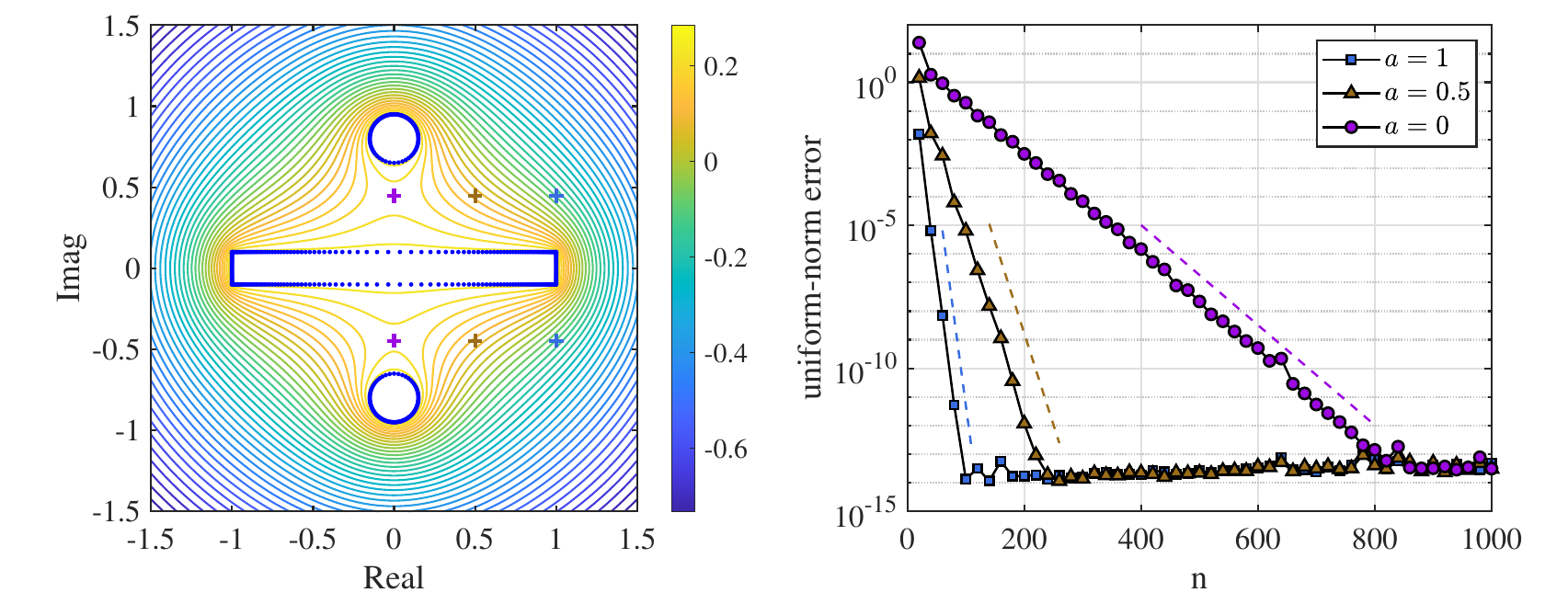}
  \caption{Left: Interpolated nodes (blue) and contours of the discrete potential $U_{\mu_n}$ in the "÷" shaped region when $n=500$.
  Right: The convergence rates compared with the theoretical convergence rates.}
  \label{fig:poly1}
\end{figure}

The polynomial interpolation of BIEP can effectively handle disjoint regions in the complex plane. In our second example, we consider a "÷" shaped region, depicted in the left panel of \cref{fig:poly1}. The objective function is defined as $f=1/(0.2+(x-a)^2)$, with $a=0$, $0.5$, and $1$. Similar to the first example, we can derive the theoretical rate of convergence based on the singularity locations. The right panel of \cref{fig:poly1} illustrates the convergence behavior, contrasting the actual speed with the theoretical prediction. The vertical axis represents the maximum relative error on the "÷" shaped region. As the interpolation order $n$ increases, the convergence speed progressively aligns with the theoretical expectation. Ultimately, both the actual and theoretical convergence rates achieve an accuracy of approximately $10^{-14}$ orders of magnitude.

\subsection{Rational interpolants for near-singularities functions}
\label{sec:4.2}
In this subsection, we investigate the interpolation of analytic functions with singularities near the interval $E=[-1,1]$. We compare our method with the conformal mapping technique proposed by Tee and Trefethen \cite{Tee2006}, as well as the Floater-Hormann method that uses equidistant nodes for interpolation \cite{Floater2007}. For the latter, we utilize the implementation \emph{Chebfun(f(X),`equi')} within Chebfun \cite{Driscoll2014}. The comparison is based on $\|f-r_{n}\|_{\infty}$, the maximum absolute error of the corresponding interpolants, computed over the interval $-1:0.00001:1$.

\begin{figure}[htbp]
  \centering
  \includegraphics[width=1\linewidth]{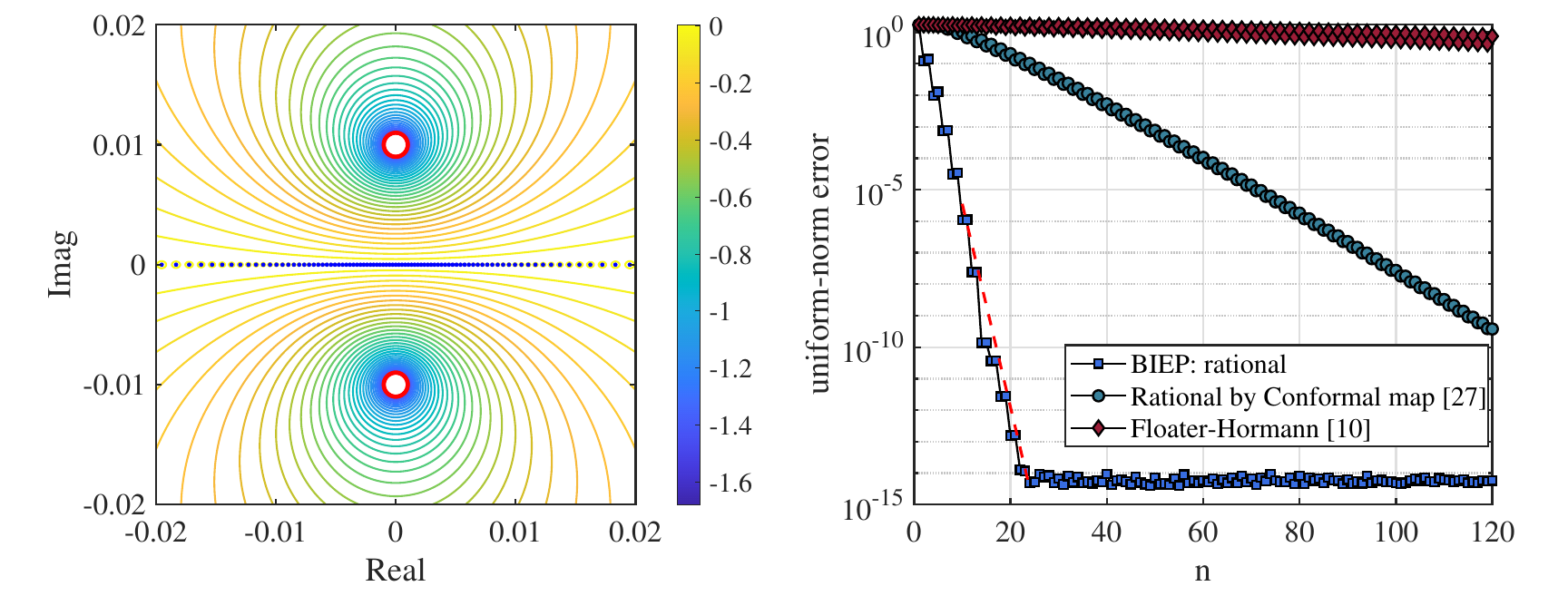}
  \caption{Left: Local potential contour with nodes and poles using our method for $n=120$.
  Right: Comparison of the effects of several types of rational interpolation on analytic functions with isolated singularities, where the red dashed line represents the theoretical rate of convergence (from \eqref{eq:2g}) of BIEP.
   }
  \label{fig:3}
\end{figure}

The third example considers $f(x) = \exp\left((1+10^4 x^2)^{-1}\right)$ with essential singularities at $\pm 0.01i$. In this case, we define the domain $F$ as two circles centered at the singularities with radius $0.001$. The distribution of poles on the boundaries of these circles encourages nodes to cluster, maintaining an equilibrium potential. This behavior is illustrated in the left panel of \cref{fig:3}. The difference in potentials $(c_1 + c_2)$ between $E$ and $F$ determines the convergence rate, as indicated by the red dashed line in the right panel. Since $f$ is analytic outside of $F$, the actual rate of convergence closely follows the exponential decay $[\exp(-c_1-c_2)]^n$.

\begin{figure}[htbp]
  \centering
  \includegraphics[width=1\linewidth]{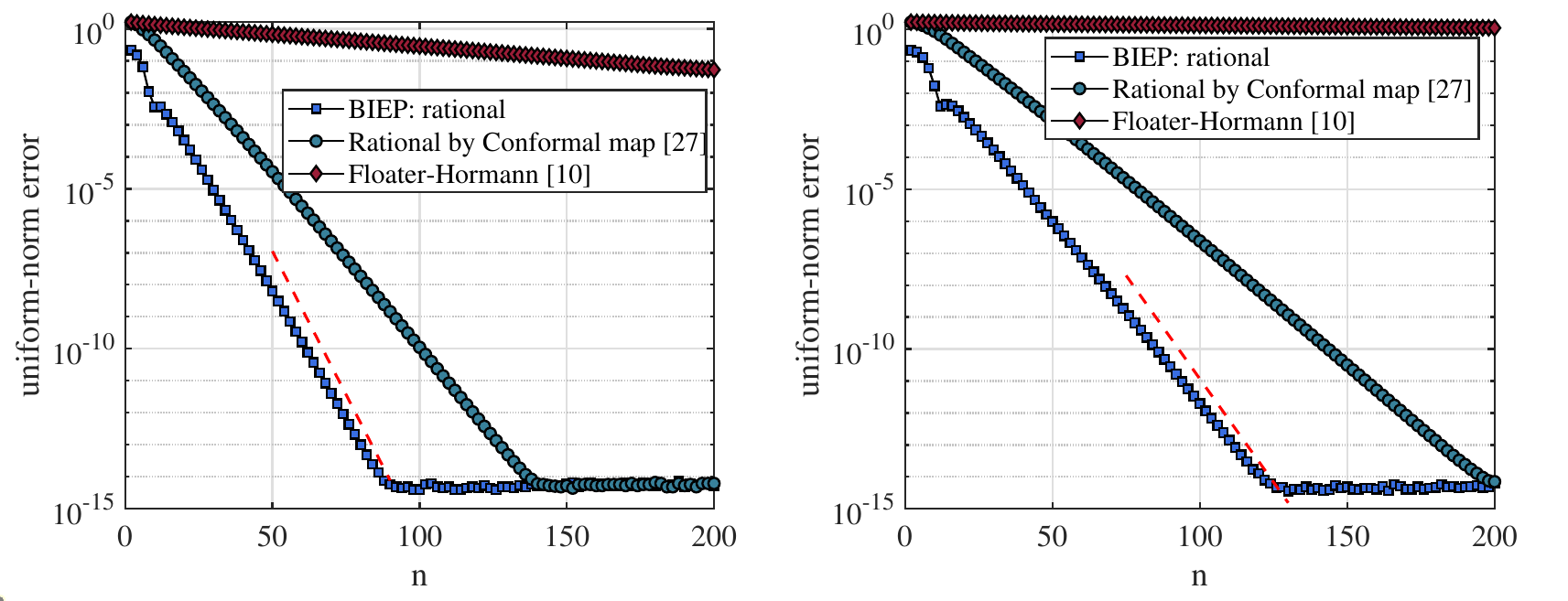}
  \caption{Uniform-norm errors for several rational interpolations of $f(x)=\exp((1+10^4 x^2)^{-0.5})$ (left) and $f(x)=\exp((1+10^6 x^2)^{-0.5})$ (right) for $n = 2:2:200$.}
  \label{fig:5}
\end{figure}

Two examples featuring branch singularities are depicted in \cref{fig:5}. The objective functions are given by $f(x)=\exp((1+10^4 x^2)^{-0.5})$ (left) and $f(x)=\exp((1+10^6 x^2)^{-0.5})$ (right). Following the second strategy outlined in \cref{sec:algorithm}, we construct $F$ using multiple line segments, each with a length of 10. Since $f$ is not analytic outside of $F$, the actual convergence speed of BIEP does not reach $[\exp(-c_1-c_2)]^n$. However, despite being much slower than in the case of an isolated singularity, this convergence speed represents an improvement of more than 50 percent compared to the speed achieved in \cite{Tee2006}.

In contrast to the polynomial interpolation of BIEP, the nodes in rational interpolation exhibit higher densities not only at the tip but also cluster in the neighborhood of $F$, as illustrated in the left panel of \cref{fig:3}. To avoid larger-scale matrices in solving the density function, one can cluster the division points in the high-density region. An adaptive mesh strategy involves first obtaining a coarse density function and then using \texttt{den2pts} with this density function to determine the division points. These points will naturally cluster in the high-density part. By calculating an appropriate approximate density function based on these division points, the value of $N$ can be effectively reduced compared to using equidistant or Chebyshev points as division points.
\section{Further discussions on the rational interpolation}\label{sec:discu}
The selection of the length parameter for the region $F$, when the cut line of the function extends to infinity, poses a challenge. In \cref{sec:3.2}, we recommended a range of 4 to 10 units if the diameter of region $E$ falls between 1 and 2 units. This broad range is justified by the fact that the actual convergence rate becomes less sensitive to this parameter as the length approaches a certain threshold.

\begin{figure}[htbp]
  \centering
  \includegraphics[width=0.55\linewidth]{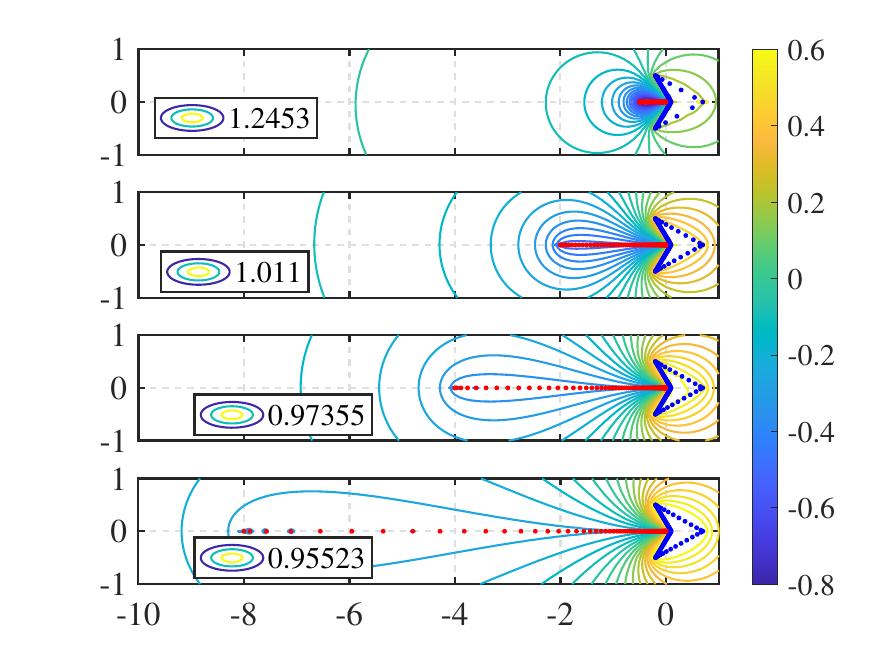}
  \caption{Potential contours of the nodes (blue) and poles (red) for $n=100$.
  The values in the figure are the potential differences $(c_1+c_2)$ corresponding to different $F$.
  }
  \label{fig:12}
\end{figure}

Here, we illustrate this with an example of interpolating the Markov function $f(z) = 1/\sqrt{z}$ on the quadrilateral $[0.1, -0.2 + 0.5i, 0.7, -0.2 - 0.5i]$. \cref{fig:12} shows the potential maps when the length of $F$ is 0.5, 2, 4, and 8, from top to bottom.

\begin{figure}[htbp]
  \centering
  \includegraphics[width=1\linewidth]{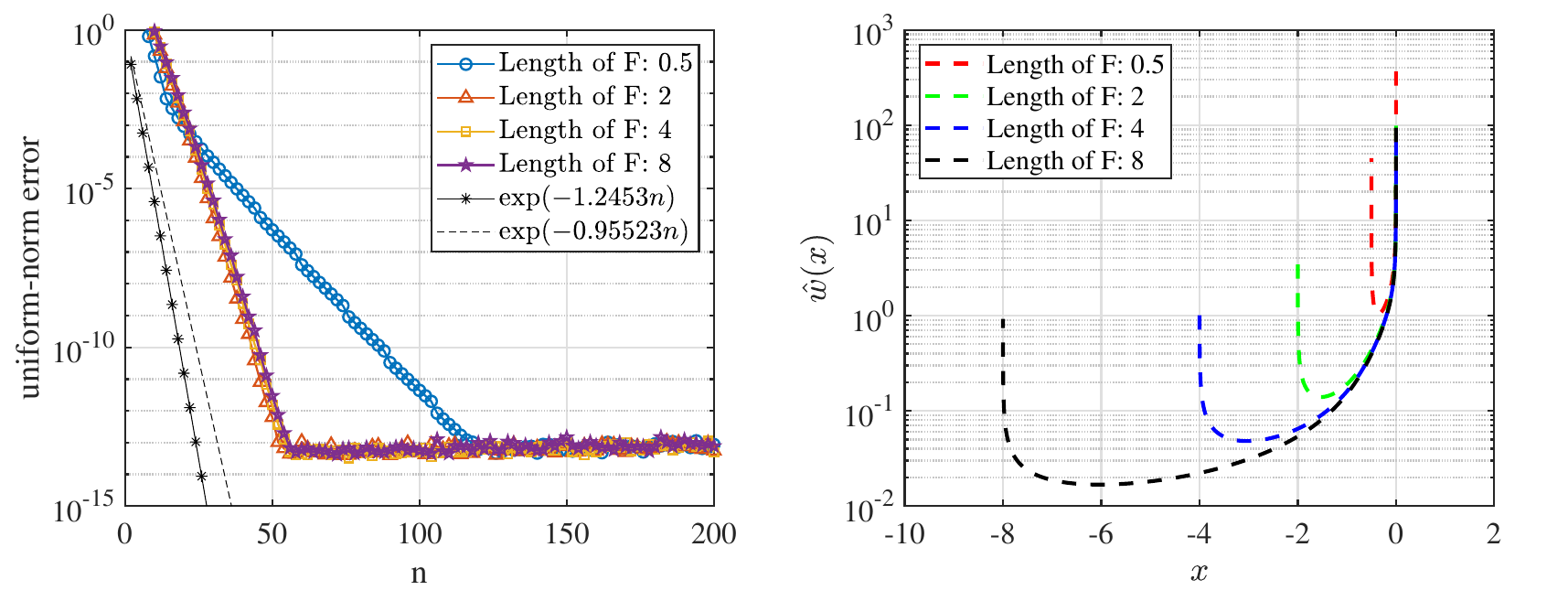}
  \caption{Left: The actual convergence speed of BIEP for different lengths of $F$ corresponding to \cref{fig:12}. 
  Right: The approximate density function $\hat{w}_F$ on different $F$.}
  \label{fig:13}
\end{figure}

It is worth noting that the shorter the length of $F$ and the more concentrated the poles are, the larger the potential difference between $E$ and $F$ will be. However, since the objective function is not analytic outside $F$, the actual convergence speed does not reach $[\exp{(-c_1-c_2)}]^n$. For this case, the actual convergence speed should be 
\begin{equation*}
\limsup_{n\to\infty} \left(\|r_n-f\|_E\right)^{\frac{1}{n}} \le \exp(\min_{\Gamma}\max_{z\in\Gamma}U(z)-c_1),
\end{equation*}
where $f$ is analytic inside $\Gamma$. The longer $F$ is, the closer the actual convergence speed is 
to $[\exp{(-c_1-c_2)}]^n$. However, $(c_1+c_2)$ becomes smaller as the length grows. The left panel of \cref{fig:13} 
illustrates the actual rate of convergence for different lengths of $F$. The errors are very close except for the 
case where the length of $F$ is 0.5.

The right panel of \cref{fig:13} illustrates the values of the approximate density function on $F$. 
The density functions behave similarly near 0, except for the case of 0.5, where the density at the distant end differs 
significantly from the density near 0 in magnitude. When the approximate density function $\hat{w}_F$ differs by 
at least one order of magnitude between its ends, the convergence rate is often acceptable. Hence, the use of longer $F$ 
is generally beneficial. We typically use lengths of 4-10 units when the diameter of $E$ is 1-2 units to handle 
such situations. Lengths shorter than 4 units may affect the speed of convergence, while lengths longer than 10 units 
are unnecessary.

We consider an example with both isolated singularities and branch points. The boundary curve
of the interpolation region $E$ is 
\[
  z(\theta)=(0.6+0.3\cos(4\theta+\pi))\exp(i\theta),\,\, (0\leq\theta\leq\pi).
\]
The function 
\(
f(z)=\sqrt{z+0.5}/((z^2+0.25)(z-0.5))
\)
has a branch cut $[-\infty,-0.5]$ and three isolated singularities.
We use $[-1,-0.5]$, $[-4.5,-0.5]$ and $[-8.5,-0.5]$ for the branch cut, respectively. 
The three isolated singularities are covered by $D(0.5i,0.01)\cup D(0.5,0.01)\cup D(-0.5i,0.01)$ in all cases.

\begin{figure}[hbpt]
\centering
\includegraphics[width=1\linewidth]{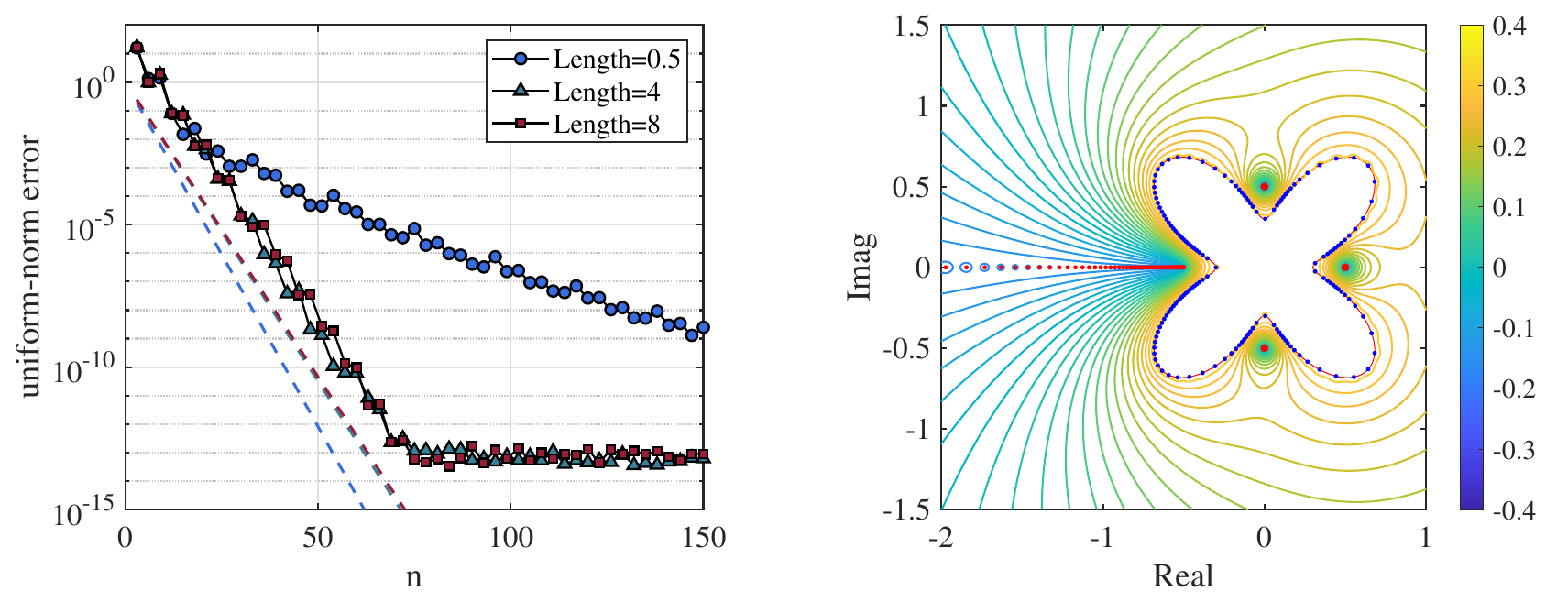}
\caption{Left: Actual convergence rate of BIEP and $[\exp(-c_1-c_2)]^n$ (dashed line) for different lengths of $F$. 
Right: Nodes (blue), poles (red), and contours of the discrete potential $U_{\mu_n}$ when the length of $F$ is 8 and $n=150$.
 }
 \label{fig:15}
\end{figure}

\Cref{fig:15} illustrates the effect of BIEP and the equipotential distribution of poles and nodes. 
For this mixed case, the actual approximations are similar and close to the convergence rate of 
$[\exp{(-c_1-c_2)}]^n$ when the length of $F$ is $4$ and $8$. The phenomenon is consistent with the previous example. 
This example also shows that BIEP is still valid for non-circular curved boundaries.

\section{Application: Solving the Laplace equation}\label{sec:Laplace}
In this section, we introduce the application of BIEP, which is characterized by its wide applicability, especially in the bounded region of the complex plane. Therefore, many applications related to rational interpolation can be realized by BIEP. We will show one of the simple applications, namely, solving 2-D Laplace boundary-value problems: 
\[
\begin{cases}
    \Delta u(z)=0, & z\in\Omega, \\
    u(z) = h(z), & z\in\partial\Omega,
\end{cases}
\]
where $\partial\Omega$ is a piecewise smooth Jordan curve.

Laplace equations, being the simplest class of elliptic partial differential equations, have been the focus of numerous numerical methods \cite{Bremer2010,Stefano2023,Gopal20192,Hochman2013,Hoskins2019}. In the two-dimensional case, a succinct approach is to convert the problem into one of approximating a harmonic function~\cite{Gopal20192}. Specifically, one can use polynomials or rationals to approximate an analytic function $f$ on the region $\Omega$, where the real part of $f$ is denoted by $u$. Subsequently, $u$ can be expressed in terms of the real part of $r$.

Approximating the analytic function $f$ in the complex plane using BIEP requires knowing the function values $f_i$ at the nodes $x_i$. 
The real part of $f_i$ is provided by the boundary conditions, leaving only the imaginary part of $f_i$ unknown. 
Therefore, to obtain the approximation for $u$, we only need the information about the imaginary part of $f_i$, 
which can be expressed as a linear combination of the Lagrange basis functions $I_i$, i.e., 
\[
  f(x)\approx r_n(x),\quad r_n:=\sum_{i=0}^{n}\left(\text{Im}(f_i)+h(x_i)\right)I_i, \quad x\in \Omega.
\]

Let $\{y_j\}_{j=0}^{K}$ be the set of points on the boundary, and $\text{Im}(f_i)$ should be such that $Re(r_n)$ 
approximates $h$ on $\{y_j\}_{j=0}^{K}$. Thus, it is only necessary to add a least squares step to the BIEP 
to solve 
\begin{equation*}
  \text{Im}\left( \begin{array}{ccc}
    I_0(y_0) & \cdots & I_n(y_0)\\
    \vdots   & \ddots & \vdots  \\
    I_0(y_K) & \cdots & I_n(y_K)\\
  \end{array} \right)
  \left( \begin{array}{c}
    \text{Im}(f_0)\\
    \vdots \\ 
    \text{Im}(f_n)\\
  \end{array} \right) =
  \text{Re}\left( \begin{array}{c}
    \sum_{i=0}^{n}h(x_i)I_i(y_0)-h(y_0)\\
    \vdots \\ 
    \sum_{i=0}^{n}h(x_i)I_i(y_K)-h(y_K)\\
  \end{array} \right).
\end{equation*}

\begin{figure}[hbpt]
\centering
\includegraphics[width=1\linewidth]{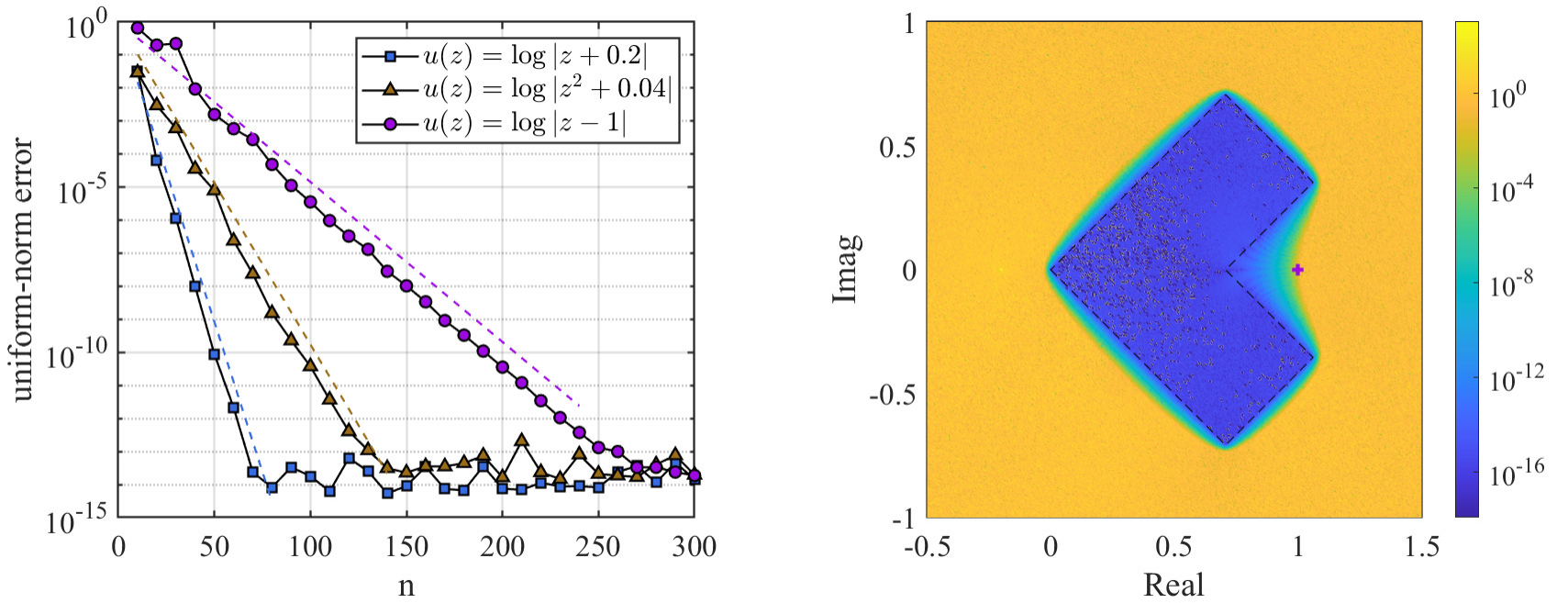}
\caption{Left: Uniform norm error for approximating several harmonic functions using the real part of polynomial interpolation, for different $n$. Right: Pointwise error for $u(z)=\log|z-1|$ when $n=300$. The violet ``+" marks the location of the singularity.
}
 \label{Figure_11}
\end{figure}

Based on the first polynomial interpolation example (\cref{fig:poly}), \cref{Figure_11} shows how 
well the real part $\text{Re}(p_n)$ of the polynomial interpolated $p_n$ approximates several harmonic functions 
with the addition of the least squares step described above. The singularities of these several harmonic functions 
are located at the same positions as the objective function in Example 1. Thus, the left panel of \cref{Figure_11} 
shows a convergence rate similar to the theoretical result \eqref{eq:1.1e} for analytic functions.

The right panel of \cref{Figure_11} illustrates the pointwise error image of the real part of the polynomial interpolation for $u(z) = \log|z-1|$ when $n = 300$. Since the maximum error occurs on the boundary, the accuracy inside the region is manageable. In this example, the maximum error appears on the boundary near the singularity.

Convergence using polynomial interpolation can be slow for data bounded near singularities or on smooth boundaries with small curvature, due to the effect of near-singularities. For these cases, better results can be achieved using rational interpolation. However, the BIEP requires prior knowledge of the singularity type and location in order to set up the F-region. Therefore, it needs to be paired with an efficient algorithm for singularity localization, such as the AAA rational method \cite{Nakatsukasa2018}.

\begin{figure}[hbpt]
\centering
\includegraphics[width=1\linewidth]{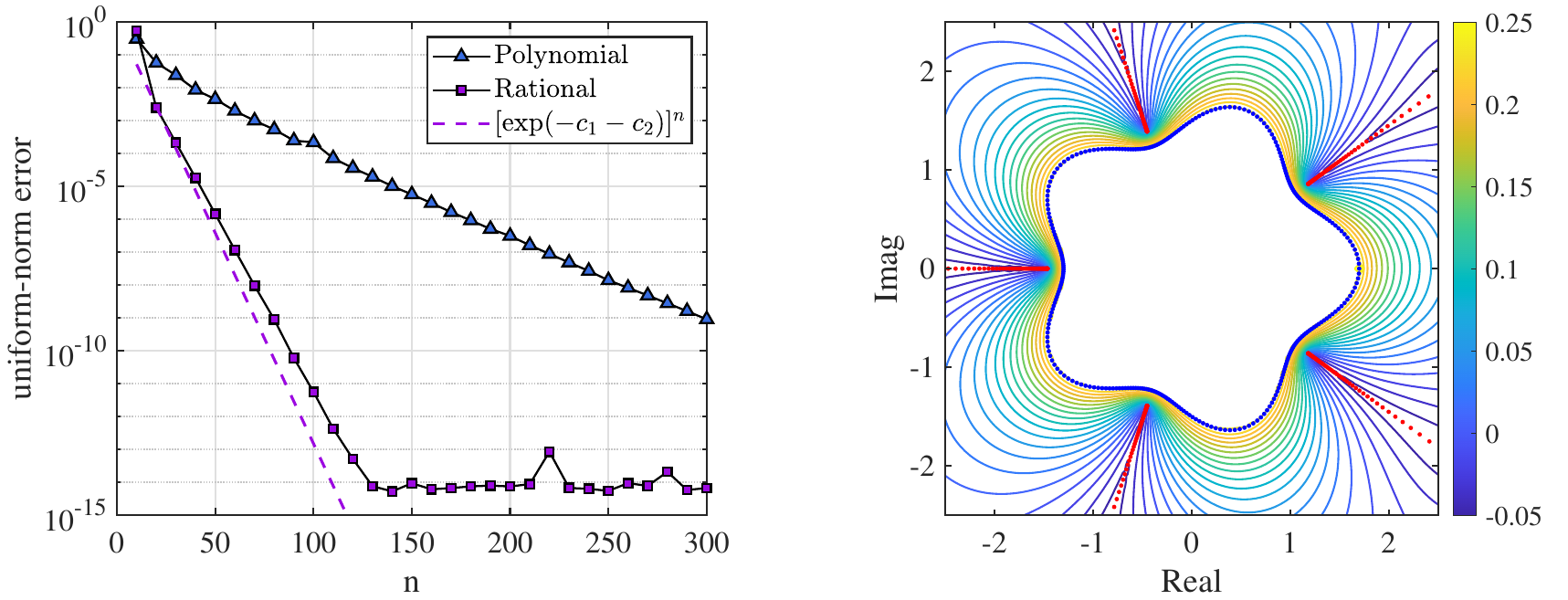}
\caption{Left: The uniform-norm error for the Laplace boundary-value problem on a smooth boundary. This error is computed using the real part of the polynomial or rational interpolation and the boundary data. Right: The nodes (blue) and poles (red) given by the rational BIEP when $n=300$, along with the contours of the discrete potential $U_{\mu_n}$.
 }
 \label{Figure_12}
\end{figure}

Here we consider an example of an analytic boundary. The boundary curve is given by 
\[z(\theta)=1.5+0.2\cos(5\theta)\exp(\mathrm{i}\theta),\quad \theta\in[0,2\pi]\]
and the boundary conditions are given by \(h(z) = \sin(3\, \mathrm{Im}(\log{z}))\).
When using the real part of the polynomial, the speed of convergence of the error is shown in the left panel of \cref{Figure_12}. 
The convergence rate of the polynomial is limited due to the ``crowding phenomenon" \cite{Gopal2019,xue2023}. However, rational interpolation 
using BIEP achieves faster convergence. Similar to the branch singularity example shown in \cref{sec:4.2}, 
the actual convergence of the error is slightly slower than \([\exp(-c_1-c_2)]^n\). Specifically, after locating the 
singularity position using the AAA method, we set up the F-region according to Strategy 2 in \cref{sec:3.2}. 
The right panel of \cref{Figure_12} illustrates the distribution of nodes and poles when \(n = 300\).

\section{Conclusion}\label{sec:conclusion}
This paper introduces efficient algorithms for polynomial and rational interpolation in the complex plane region bounded by piecewise smooth Jordan curves. The nodes and poles of these interpolation methods closely follow the equipotential distribution, and the approximation results for analytic functions align with established theoretical convergence rates. Rational interpolation, in particular, demonstrates a rapid convergence rate for nearly singular functions. Leveraging these characteristics, we propose a straightforward computational approach for solving the 2-D Laplace boundary-value problem.

For 2-D Laplace boundary-value problems, the occurrence of branching points on the boundary presents a particularly interesting case. However, the BIEP method proposed in this paper is limited to handling analytic functions and near-singular functions. Therefore, one of our next objectives is to extend the BIEP method to accommodate analytic functions with bounded weak singularities.

The core of the algorithm proposed in this paper is to generate discrete points in a complex plane region such that the discrete potential approximates the logarithmic equilibrium potential. These points have applications not only in interpolation but also in numerical algebra. The approach can be extended to solve other problems, such as the Zolotarev problem in the complex plane or to determine the optimal parameters for the alternating direction implicit (ADI) method \cite{Istace1995, Bailly2000}. Further research will explore these avenues in future work.

\section*{Code availability}
The code required to replicate the numerical examples can be found on GitHub at \href{https://github.com/Clone-Z/BIEP_example}{https://github.com/Clone-Z/BIEP\_example}.

\section*{Appendix A}\label{secA}
{\bf Proof of \cref{thm:den}}: We will first prove \eqref{eq:pdis}. Assume that  for any positive integer $K$,
there exists $n>K$ such that $\max_{1\leq i\leq n}\lvert x_{i-1}^{(n)}-x_i^{(n)}\rvert>2d$ for some positive constant $d$. Divide $\partial E$ into
$m=\left\lceil \lvert\partial E\rvert _L/d \right\rceil$ parts $\Gamma_1,\Gamma_2,\cdots,\Gamma_m$ with equal length, where $\left\lceil z\right\rceil$ denotes the smallest integer greater than or equal to $z$.
Then, for any positive integer $K$, there exists $n(>K)$ and some $i(n)$ with $1\leq i(n)\leq m$ such that
$n_{\Gamma_{i(n)}}=0$. Since $\{\Gamma_i\}_{i=1}^m$ are a finite set, thus there is $i_0$ with $1\leq i_0\leq m$ and a subsequence such that $n_{\Gamma_{i(n_k)}}=n_{\Gamma_{i_0}}=0$ for $k=1,2,\ldots$, which implies that (\ref{eq:3a}) does not hold, i.e., $\{x_i^{(n)}\}_{i=0}^n$ does not
obey $w(t)$ on $\partial E$ since $\lim_{k\rightarrow \infty}\frac{n_{\Gamma_{i(n_k)}}}{n_k}=\int_{\Gamma_{i_0}}w(t)\,\lvert\mathrm{d}t\rvert=\lim_{k\rightarrow \infty}\frac{n_{\Gamma_{i_0}}}{n_k}=0$, which is contradicted with $\int_{\Gamma_{i_0}}w(t)\,\lvert\mathrm{d}t\rvert>0$. Then the identity \eqref{eq:pdis} holds.

Next we will show that
\begin{equation}
  \lim_{n\to\infty}\int g \, \mathrm{d}\mu_{n}=\int g \, \mathrm{d}\mu
\end{equation}
for arbitrary continuous function $g$ on $\partial E$, where $\mu_{n} = 1/(n+1)\sum_{i=0}^{n}\delta_{x_i^{(n)}}$.

Suppose $E$ is a union of  simply or multicoonencted sets and $\partial E = \bigcup_{i\in I} \partial E^{i}$, where $\partial E^i$ is close and smooth
and $I$ is a finite index set. Then we just need to prove that weak* convergence holds on every $\partial E^i$.

Denote $\max_{t\in\partial E^i}\lvert g(t)\rvert=G_i$.
Let the parametric curve of the smooth curve $\partial E^i$ about the arc length $s$ be $T_i(s)$
($0\leq s\leq S_i=\lvert\partial E^i\rvert _L$). Then we have
\[
  \int_{\partial E^{i}} g(t)\, \mathrm{d}\mu(t) = \int_{\partial E^{i}} g(t)w(t)\,\lvert\mathrm{d}t\rvert
  = \int_{0}^{S_i} g(T_i(s))w(T_i(s))\, \mathrm{d}s.
\]
Let $0=s_0<s_1<\cdots<s_m=S_i$ satisfy
\begin{equation} \label{eq:A3}
  \int_{s_{k-1}}^{s_{k}} w(T_i(s)) ds=\frac{W_i}{m},\quad W_i=\int_{0}^{S_i} w(T_i(s))\, \mathrm{d}s\leq 1 .
\end{equation}
Since $w>0$, we have $\lim_{m\to\infty}\max_{1\leq k\leq m}\lvert s_k-s_{k-1}\rvert=0$ from Identity  \eqref{eq:pdis}.

Note that $g(T_i(s))$ is continuous for $s\in[0,S_i]$. Then for arbitrary $\varepsilon>0$, there exists $K$ such that
\begin{equation}\label{eq:A4}
  \max_{s\in[s_{k-1},s_k]}g(T_i(s))-\min_{s\in[s_{k-1},s_k]}g(T_i(s))<\frac{\varepsilon}{2},\quad 1\leq k\leq m,\quad m\geq K.
\end{equation}
Given an $m$  ($\geq M$), by applying the mean value theorem of integration, there exists $\xi_k\in[s_{k-1},s_k]$ such that
\begin{equation} \label{eq:A5}
  \int_{0}^{S_i} g(T_i(s))w(T_i(s))\, \mathrm{d}s = \frac{W_i}{m}\sum_{k=1}^{m} g(T_i(\xi_k)).
\end{equation}
Furthermore, from
\[
  \int_{s_{k-1}}^{s_k} g(T_i(s))\, \mathrm{d}\mu_{n+1}(T_i(s))=\frac{\sum_{x\in X_i} g(x)}{n+1}
\]
where $X_i=\{x_i^{(n)}\}_{i=0}^n\cap T_i([s_{k-1},s_k])$,
there  exists $\zeta_k\in[s_{k-1},s_k]$ such that
\[
  \frac{\sum_{x\in X_i} g(x)}{n+1}=g(\zeta_k)\frac{n_{\partial E[\widehat{t_{k-1}^i,t_{k}^i}]}}{n+1}
\]
and
\begin{equation} \label{eq:A8}
  \lvert g(T_i(\xi_k))-g(T_i(\zeta_k))\rvert <\frac{\varepsilon}{2},
\end{equation}
where $t_{k-1}^i=T_i(s_{k-1})$ and $t_{k}^i=T_i(s_{k})$.

In addition, from (\ref{eq:3a}), there  exists a positive integer $N$ such that for $n>N$
\begin{equation} \label{eq:A6}
  \left\lvert\frac{n_{\partial E[\widehat{t_{k-1}^i,t_{k}^i}]}}{n+1}-\int_{s_{k-1}}^{s_{k}} w(T_i(s))\, \mathrm{d}s \right\rvert
  \leq \frac{\varepsilon m}{2G_i W_i},\, k=1,2,\ldots,m,
\end{equation}
which, together with (\ref{eq:A3}) and (\ref{eq:A6}), deduces
\begin{equation} \label{eq:A7}
  \frac{n_{\partial E[\widehat{t_{k-1}^i,t_{k}^i}]}}{n+1}\Big/\int_{s_{k-1}}^{s_{k}} w(T_i(s))\, \mathrm{d}s
  \in [1-\frac{\varepsilon}{2G_i},1+\frac{\varepsilon}{2G_i}].
\end{equation}
Thus, together with (\ref{eq:A5}), (\ref{eq:A8}) and (\ref{eq:A7}),  we have for $n>N$ that
\begin{align*}
 \left\lvert\int_{\partial E^{i}}g\, \mathrm{d}\mu_{n+1}-\int_{\partial E^{i}}g\, \mathrm{d}\mu\right\rvert
 &=\left\lvert \sum_{k=1}^{m}g(T_i(\zeta_k))\frac{n_{\partial E[\widehat{t_{k-1}^i,t_{k}^i}]}}{n+1}-\frac{W_i}{m}\sum_{k=1}^{m} g(T_i(\xi_k))\right\rvert\\
 &\leq\frac{W_i}{m}\sum_{k=1}^{m}(\left\lvert g(T_i(\xi_k))-g(T_i(\zeta_k))\right\rvert+\frac{\varepsilon}{2G_i}\left\lvert g(T_i(\zeta_k))\right\rvert)\\
 &\leq \frac{W_i}{m}\sum_{k=1}^{m}(\frac{\varepsilon}{2}+\frac{\varepsilon}{2})\\
 &\leq \varepsilon.
\end{align*}
This completes the proof.


\bibliographystyle{siamplain}
\bibliography{references}
\end{document}